\documentclass[12pt]{amsart}
\usepackage{amssymb}
\usepackage{array}
\usepackage{hyperref}

\newtheorem{theorem}{Theorem}[section]
\newtheorem{definition}[theorem]{Definition}
\newtheorem{proposition}[theorem]{Proposition}
\newtheorem{corollary}[theorem]{Corollary}

\newtheorem{conjecture}[theorem]{Conjecture}
\newtheorem{problem}[theorem]{Problem}
\newtheorem{question}[theorem]{Question}
\newtheorem{fact}[theorem]{Fact}

\begin{document}

\title[Quantum permutations and Hadamard matrices]
{Quantum permutations, Hadamard matrices, and the search for matrix models}

\author{Teodor Banica}
\address{T.B.: Department of Mathematics, Cergy-Pontoise University, 95000 Cergy-Pontoise, France. {\tt teodor.banica@u-cergy.fr}}


\begin{abstract}
This is a presentation of recent work on quantum permutation groups, complex Hadamard matrices, and the connections between them. A long list of problems is included. We include as well some conjectural statements, about matrix models.
\end{abstract}

\maketitle

{\em Dedicated to Professor S.L. Woronowicz, on the occasion of his 70-th birthday.}

\tableofcontents

\section*{Introduction}

A complex Hadamard matrix is a square matrix $H\in M_n(\mathbb C)$ whose entries are on the unit circle, $|H_{ij}|=1$, and whose rows are pairwise orthogonal, when regarded as elements of $\mathbb C^n$. It follows from definitions that the columns are pairwise orthogonal as well.

The basic example is the Fourier matrix $F_n=(w^{ij})$, where $w=e^{2\pi i/n}$. The name comes from the fact that $F_n$ is the matrix of the discrete Fourier transform, over $\mathbb Z_n$. This is in fact the one and only basic example, typically the only known one when $n$ is prime.

However, lots of other interesting Hadamard matrices, each one having its own story, do exist. For instance in the real case the Hadamard conjecture, going back to the 19th century, states that there is a such a matrix $H\in M_n(\mathbb R)$, for any $n\in 4\mathbb N$.

Due to their remarkable combinatorial properties, the complex Hadamard matrices appear in a wealth of concrete situations, in connection with subfactors, spin models, knot invariants, planar algebras, quantum groups, and various aspects of combinatorics, functional analysis, representation theory, and quantum physics.

The present article is a survey on complex Hadamard matrices, quantum permutation groups, and the connections between them. The idea is very simple:
\begin{enumerate}
\item First, the Fourier matrix $F_n$ is by definition cyclic, and is in fact naturally associated to $\mathbb Z_n$. This suggests that any complex Hadamard matrix $H$ should be in fact associated to a group-like object $G$, describing its symmetries.

\item A glimpse at the various complex Hadamard matrices (different from $F_n$) which have been constructed leads to the quick conclusion that their symmetries are rather of ``alien'' nature, and cannot be described by a usual group $G$.

\item So, we should rather look for a quantum group. And, with this idea in mind, the solution appears: the orthogonality condition on the rows of $H$ reminds the ``magic'' condition on the coordinates of the quantum permutation group $S_n^+$.

\item With a bit more work, one can associate to any complex Hadamard matrix $H\in M_n(\mathbb C)$ a certain quantum permutation group $G\subset S_n^+$. As a first example, associated in this way to the Fourier matrix $F_n$ is the usual cyclic group $\mathbb Z_n$.
\end{enumerate}

There's perhaps at this point some need for explanations: the complex Hadamard matrices are such concrete and beautiful objects, so why doing all this?

The problem is that any complex Hadamard matrix $H$ produces a certain sequence of numbers $c_k\in\mathbb N$, which are extremely hard to compute. These numbers, called ``quantum invariants'', appeared from the work of Popa and Jones on subfactor theory, and play a key role in most of the quantum algebraic contexts where $H$ appears. Their formal definition is $c_k=\dim(P_k)$, where $P=(P_k)$ is the planar algebra associated to $H$.

The point now is that, with the above approach, the numbers $c_k$ have a very simple interpretation. Indeed, they appear as moments of the main character of $G$:
$$c_k=\int_GTr(g)^k\,dg$$

This formula is of course quite encouraging. For instance it shows that, from the point of view of quantum algebra in a very large sense, $G$ is indeed the correct symmetry group of $H$. Of course, these results are quite new, and the correspondence $H\to G$ still needs to be developed in order to reach to some concrete results. It is our hope that so will be the case, and in the present article we explain the state-of-art of the subject.

The paper is organized as follows: in 1-3 we discuss the quantum permutation groups, in 4 we discuss the complex Hadamard matrices, in 5 we describe the correspondence and its basic properties, and in 6 we present some random matrix speculations.

\subsection*{Acknowledgements}

It would be impossible to make a full list of those who have made this theory possible, but at the beginning we certainly have Professor S.L. Woronowicz, to whom this article is dedicated. It is also a pleasure to thank Julien Bichon, Beno\^it Collins and Jean-Marc Schlenker for their key contribution to the subject, Stephen Curran, Sonia Natale, Adam Skalski and Roland Speicher for their work on quantum permutations, and Dietmar Bisch, Vaughan Jones and Remus Nicoara for various subfactor discussions. 

The present article is based on lecture notes from minicourses given at the University of Caen in Fall 2010, and at the Erwin Schr\"odinger Institute in Vienna in Spring 2011, and I am grateful to Roland Vergnioux and Dan Voiculescu for the invitations.

This work was supported by the ANR grants ``Galosint'' and ``Granma''.

\section{Quantum permutations}

Consider the embedding $S_n\subset O_n$ given by the usual permutation matrices, $\sigma(e_j)=e_{\sigma(j)}$. The matrix entries of a permutation $\sigma\in S_n$ are then:
$$\sigma_{ij}
=\begin{cases}
1&{\rm if}\ \sigma(j)=i\\
0&{\rm if\  not}
\end{cases}$$ 

Let us introduce now the functions $u_{ij}:S_n\to\mathbb C$ given by $u_{ij}(\sigma)=\sigma_{ij}$. These are the $n^2$ coordinates on $S_n$, viewed as an algebraic group. 

We can see that the matrix $u=(u_{ij})$ is magic, in the following sense:

\begin{definition}
A square matrix $u\in M_n(A)$ is called ``magic'' if all its entries are projections, summing up to $1$ on each row and column.
\end{definition}

Here of course $A$ is a $C^*$-algebra, and a projection is by definition an element $p\in A$ satisfying $p=p^*=p^2$. In the above example, the algebra is $A=C(S_n)$.

Consider now the multiplication map $(\sigma,\tau)\to\sigma\tau$. At the level of functions, this corresponds to a map $\Delta:A\to A\otimes A$ called comultiplication, given by $\Delta f(\sigma,\tau)=f(\sigma\tau)$. On the above elements $u_{ij}$, the action of $\Delta$ is extremely simple:
$$\Delta(u_{ij})=\sum_{k=1}^nu_{ik}\otimes u_{kj}$$

In addition, we know from the Stone-Weierstrass theorem that the coordinate functions  $u_{ij}$ generate $A$. With a bit more work one can prove that these generators $u_{ij}$ satisfy no other relations besides the magic ones, so we have reached to the following conclusion:

\begin{theorem}
$C(S_n)$ is the universal commutative $C^*$-algebra generated by the entries of a $n\times n$ magic matrix $u$, and its comultiplication is given by $\Delta(u_{ij})=\sum_ku_{ik}\otimes u_{kj}$.
\end{theorem}

Summarizing, we have a fully satisfactory description of the Hopf algebra $C(S_n)$. Observe that the counit map $\varepsilon(f)=f(1)$ and the antipode map $Sf(\sigma)=f(\sigma^{-1})$ are given as well by very simple formulae on the standard generators, namely:
$$\varepsilon(u_{ij})=\delta_{ij},\quad\quad S(u_{ij})=u_{ji}$$

We recall now that a compact quantum group is an abstract object $G$, having no points in general, but which is described by a well-defined Hopf $C^*$-algebra $A=C(G)$. 

The axioms for Hopf $C^*$-algebras, found by Woronowicz in \cite{wo3}, are as follows:

\begin{definition}
A Hopf $C^*$-algebra is a $C^*$-algebra $A$, given with a morphism of $C^*$-algebras $\Delta:A\to A\otimes A$, called comultiplication, subject to the following  conditions:
\begin{enumerate}
\item Coassociativity: $(\Delta\otimes id)\Delta=(id\otimes\Delta)\Delta$.

\item $\overline{span}\,\Delta(A)(A\otimes 1)=\overline{span}\,\Delta(A)(1\otimes A)=A\otimes A$.
\end{enumerate}
\end{definition}

The basic example is $A=C(G)$, where $G$ is a compact group, with $\Delta f(g,h)=f(gh)$. The fact that $\Delta$ is coassociative corresponds to $(gh)k=g(hk)$, and the conditions in (2) correspond to the cancellation rules $gh=gk\implies h=k$ and $gh=kh\implies g=k$. 

Conversely, any commutative Hopf $C^*$-algebra is of the form $C(G)$. Indeed, by the Gelfand theorem we have $A=C(G)$, with $G$ compact space, and (1,2) above tell us that $G$ is a semigroup with cancellation. By a well-known result, $G$ follows to be a group.

The following key definition is due to Wang \cite{wa2}:

\begin{definition}
$C(S_n^+)$ is the universal $C^*$-algebra generated by the entries of a $n\times n$ magic matrix $u$, with comultiplication $\Delta(u_{ij})=\sum_ku_{ik}\otimes u_{kj}$.
\end{definition}

Observe the similarity with Theorem 1.2. Wang proved that this algebra satisfies Woronowicz's axioms in Definition 1.3, so that $S_n^+$ is a compact quantum group.

The counit and antipode are given by the same formulae as in the classical case, namely $\varepsilon(u_{ij})=\delta_{ij}$ and $S(u_{ij})=u_{ji}$. Observe that the square of the antipode is the identity:
$$S^2=id$$ 

As a first remark, we have a surjective morphism $C(S_n^+)\to C(S_n)$, which corresponds to a quantum group embedding $S_n\subset S_n^+$. The very first problem is whether this embedding is an isomorphism or not, and the answer here is as follows:

\begin{proposition}
The quantum groups $S_n^+$ are as follows:
\begin{enumerate}
\item At $n=1,2,3$ we have $S_n=S_n^+$.

\item At $n\geq 4$, $S_n^+$ is a non-classical, infinite compact quantum group.
\end{enumerate}
\end{proposition}

This is clear indeed at $n=2$, where each magic matrix must be of the form:
$$u=\begin{pmatrix}p&1-p\\1-p&p\end{pmatrix}$$

At $n=3$ the statement is that any 9 projections forming a $3\times 3$ magic matrix commute. This is not obvious, but can be seen with a Fourier over $\mathbb Z_3$ base change.

At $n=4$ now, consider the following magic matrix:
$$u=\begin{pmatrix}p&1-p&0&0\\1-p&p&0&0\\0&0&q&1-q\\0&0&1-q&q\end{pmatrix}$$

Here $p,q$ can be any projections. If we choose $p,q\in B(H)$ to be free, then $<p,q>$ is not commutative and infinite dimensional, and so must be $C(S_4^+)$.

Finally, at $n\geq 5$ the result follows fom $S_k^+\subset S_n^+$ for $k\leq n$.

The following result was proved as well by Wang in \cite{wa2}:

\begin{theorem}
$C(S_n^+)$ is the universal Hopf $C^*$-algebra coacting on $\mathbb C^n$. In other words, $S_n^+$ is the universal compact quantum group acting on $\{1,\ldots,n\}$. 
\end{theorem}

The idea here is that a coaction map $\alpha(e_i)=\sum_je_j\otimes u_{ji}$ is a morphism of algebras if and only if the matrix $u=(u_{ij})$ is magic. Thus, if we try to construct the universal algebra coacting on $\mathbb C^n$, we are naturally led to Definition 1.4. See Wang \cite{wa2}.

Let us go back now to the above proof of the fact that $S_4^+$ is infinite. If we denote by $X$ the ``quantum space'' corresponding to the algebra $<p,q>$, the above proof is:
$$X\subset S_4^+,\ |X|=\infty\implies|S_4^+|=\infty$$

The point now is that $X$ is in fact a quantum group. Indeed, it is well-known that the algebra $<p,q>$ generated by two free projections is always the same, regardless of the projections $p,q$ that we choose. One model is provided by $p=(1+g)/2$ and $q=(1+h)/2$, where $g,h$ are the standard generators of $D_\infty=\mathbb Z_2*\mathbb Z_2$, and we get:

\begin{proposition}
We have an embedding $\widehat{D}_\infty\subset S_4^+$.
\end{proposition}

In order to explain the meaning of this statement, we have to go back to Definition 1.3. Recall that the basic example there was $A=C(G)$, with $G$ compact group. The point is that there is another class of ``basic examples'', namely the group algebras $A=C^*(\Gamma)$, with $\Delta(g)=g\otimes g$. So, the abstract dual $\widehat{\Gamma}$ of a discrete group $\Gamma$ can be defined by: 
$$C^*(\Gamma)=C(\widehat{\Gamma})$$

With this notion in hand, the above result simply tells us that the morphism $C(S_4^+)\to C^*(D_\infty)$ commutes with the comultiplications. But this is clear from definitions.

Proposition 1.7 raises several questions, the first one being the classification of group dual subgroups $\widehat{\Gamma}\subset S_n^+$. The answer here, found by Bichon in \cite{bi2}, is as follows:

\begin{theorem}
For any quotient group $\mathbb Z_{n_1}*\ldots*\mathbb Z_{n_k}\to\Gamma$ we have $\widehat{\Gamma}\subset S_n^+$, where $n=n_1+\ldots+n_k$. Any group dual subgroup of $S_n^+$ appears in this way.
\end{theorem}

The proof of the first assertion is similar to that of Proposition 1.7, by using a diagonal concatenation of magic matrices. The proof of the converse is quite technical, see \cite{bi2}.

Recall now Cayley's theorem, stating that any finite quantum group $G$ is a permutation group: $G\subset S_n$, with $n=|G|$. The above result, when applied to the quotient map $\mathbb Z_n^{*n}\to G$, shows that $\widehat{G}$ is a quantum permutation group: $\widehat{G}\subset S_{n^2}^+$. So, the following problem appears: is any finite quantum group a quantum permutation group?

This question was recently answered in the negative \cite{fin}:

\begin{theorem}
There exists a quantum group of finite order, namely $\mathbb Z_4\#\widehat{S}_3$ having order $24$, which is not a quantum permutation group.
\end{theorem}

We refer to \cite{fin} for the proof, as well as for a number of positive results: for instance any quantum group whose order divides $p^3$ or $pqr$ is a quantum permutation one.

The other question raised by Proposition 1.7 is that of understanding the structure of $S_4^+$. The answer here, quite surprising, involves the following definition, from \cite{qfp}:

\begin{definition}
$C(SO_3^{-1})$ is the universal $C^*$-algebra generated by the entries of a $3\times 3$ orthogonal matrix $a=(a_{ij})$, with the following relations:
\begin{enumerate}
\item Skew-commutation: $a_{ij}a_{kl}=\pm a_{kl}a_{ij}$, with sign $+$ if $i\neq k,j\neq l$, and $-$ if not.

\item Twisted determinant condition: $\Sigma_{\sigma\in S_3}a_{1\sigma(1)}a_{2\sigma(2)}a_{3\sigma(3)}=1$.
\end{enumerate}
\end{definition}

Consider indeed the matrix $a^+=diag(1,a)$, corresponding to the action of $SO_3^{-1}$ on $\mathbb C^4$, and apply to it the Fourier transform over the Klein group $K=\mathbb Z_2\times\mathbb Z_2$: 
$$u=
\frac{1}{4}
\begin{pmatrix}
1&1&1&1\\
1&-1&-1&1\\
1&-1&1&-1\\
1&1&-1&-1
\end{pmatrix}
\begin{pmatrix}
1&0&0&0\\
0&a_{11}&a_{12}&a_{13}\\
0&a_{21}&a_{22}&a_{23}\\
0&a_{31}&a_{32}&a_{33}
\end{pmatrix}
\begin{pmatrix}
1&1&1&1\\
1&-1&-1&1\\
1&-1&1&-1\\
1&1&-1&-1
\end{pmatrix}$$

The point is that this matrix is magic, and vice versa, i.e. the Fourier transform over $K$ converts the relations in Definition 1.10 into the relations in Definition 1.1. 

This gives the following result, proved in \cite{qfp}:

\begin{theorem}
We have $S_4^+=SO_3^{-1}$, with the algebraic isomorphism given on the standard generators by the Fourier transform over the Klein group $K=\mathbb Z_2\times\mathbb Z_2$.
\end{theorem}

This isomorphism raises the problem of classifying the ``quantum groups acting on 4 points'', i.e. the subgroups of $S_4^+$. The result here, obtained in \cite{qfp}, is as follows:

\begin{theorem}
The following is the list of subgroups of $S_4^+=SO_3^{-1}$:
\begin{enumerate}
\item Infinite quantum groups: $S_4^+$, $O_2^{-1}$, $\widehat{D}_\infty$.

\item Finite groups: $S_4$, and its subgroups.

\item Finite group twists: $S_4^{-1}$, $A_5^{-1}$.

\item Series of twists: $D_{2n}^{-1}$ $(n\geq 3)$, $DC^{-1}_n$ $(n\geq 2)$.

\item A group dual series: $\widehat{D}_n$, with $n\geq 3$.
\end{enumerate}
\end{theorem}

A finer version of this result is provided by the following table, jointly classifying the subgroups of $SO_3$ and of $S_4^+$, and containing as well the corresponding ADE graphs: 

\setlength{\extrarowheight}{3pt}
\begin{center}
\begin{tabular}[t]{|l|l|l|l|l|l|}
\hline {\rm Subgroup of} $SO_3$&{\rm ADE graph}&{\rm Subgroup of} $S_4^+$\\
\hline $\mathbb Z_1$, $\mathbb Z_2$, $\mathbb Z_3$, $\mathbb Z_4$&$\tilde{A}_1$, $\tilde{A}_3$, $\tilde{A}_5$, $\tilde{A}_7$&$\mathbb Z_1$, $\mathbb Z_2$, $\mathbb Z_3$, $K$\\
\hline $\mathbb Z_{2n-1}$, $\mathbb Z_{2n}$&$\tilde{A}_{4n-3}$,
$\tilde{A}_{4n-1}$&$\square$, $\widehat{D}_n$\\
\hline $SO_2$, $O_2$, $SO_3$&$A_{-\infty,\infty}$, $D_\infty$, $A_\infty$&$\widehat{D}_\infty$, $O_2^{-1}$, $S_4^+$\\
\hline $K$, $D_4$&$\tilde{D}_4$, $\tilde{D}_6$&$\mathbb Z_4$, $D_4$\\
\hline $D_{2n-1}$, $D_{2n}$&$\tilde{D}_{2n+1}$, $\tilde{D}_{2n+2}$& $\square$, $D_{2n}^{-1}/DC_n^{-1}$\\
\hline $\square$, $\square$ &$\square$, $\square$&$D_1$, $S_3$\\
\hline $A_4$, $S_4$, $A_5$&$\tilde{E}_6$, $\tilde{E}_7$, $\tilde{E}_8$&$A_4$, $S_4/S_4^{-1}$, $A_5^{-1}$\\
\hline
\end{tabular}
\end{center}
\setlength{\extrarowheight}{0pt}
\bigskip

Here the squares denote the various missing objects in the correspondence, and the slashes indicate the places where the correspondence is not one-to-one.

The meaning of the middle column is a bit formal, in the sense that the graphs there are Cayley type graphs for the corresponding groups or quantum groups. See \cite{qfp}.

\section{Representation theory}

We have seen in the previous section that $S_4^+$ is understood quite well. In order to discuss now the general case, we need some basic notions of representation theory:

\begin{definition}
A unitary representation of $G$ is described by the corresponding corepresentation of $A=C(G)$, which is a unitary matrix $u\in M_n(A)$ satisfying:
\begin{eqnarray*}
\Delta(u_{ij})&=&\sum_ku_{ik}\otimes u_{kj}\\
\varepsilon(u_{ij})&=&\delta_{ij}\\
S(u_{ij})&=&u_{ji}^*
\end{eqnarray*}
The representation is called ``fundamental'' if its coefficients generate $A$. 
\end{definition}

Woronowicz proved that the compact quantum groups are semisimple. More precisely, his results include a Peter-Weyl type theory \cite{wo1}, and a Tannaka-Krein duality \cite{wo2}.

The main problem is to classify the irreducible representations of $G$, and to find their fusion rules. In the self-adjoint case, $u=\bar{u}$, it is known that each $r\in Irr(G)$ appears in a certain tensor power $u^{\otimes k}$. So, the problem is to compute the spaces $End(u^{\otimes k})$.

It is convenient to enlarge the attention to the following family of spaces:
$$C(k,l)=Hom(u^{\otimes k},u^{\otimes l})$$

The point is that these spaces form a tensor category, and Woronowicz's Tannakian duality results in \cite{wo2} can be effectively used for their computation.

Let $P(k,l)$ be the set of partitions between $k$ upper points, and $l$ lower points.

\begin{definition}
Associated to $\pi\in P(k,l)$ is the linear map $T_\pi:(\mathbb C^n)^{\otimes k}\to(\mathbb C^n)^{\otimes l}$,
$$T(e_{i_1}\otimes\ldots\otimes e_{i_k})=\sum_{j_1\ldots j_l}\delta_\pi\begin{pmatrix}i_1&\ldots&i_k\\ j_1&\ldots&j_l\end{pmatrix}e_{j_1}\otimes\ldots\otimes e_{j_l}$$
where $\delta_\pi=1$ if the indices contained in any block of $\pi$ are all equal, and $\delta_\pi=0$ if not.
\end{definition}

It follows from definitions that the spaces $span(T_\pi|\pi\in P(k,l))$ form a tensor category. Indeed, the composition and tensor product of linear maps correspond to the vertical and horizontal concatenation of the partitions, and the $*$ operation corresponds to the upside-down turning of partitions. Moreover, it is easy to see that for any $\pi\in P(k,l)$ we have $T_\pi\in Hom(u^{\otimes k},u^{\otimes l})$, where $u$ is the fundamental representation of $S_n$.

Let also $NC(k,l)\subset P(k,l)$ be the set of noncrossing partitions.

\begin{theorem}
The spaces $Hom(u^{\otimes k},u^{\otimes l})$ for $S_n,S_n^+$ are as follows:
\begin{enumerate}
\item For $S_n$ we obtain $span(T_\pi|\pi\in P(k,l))$.

\item For $S_n^+$ we obtain $span(T_\pi|\pi\in NC(k,l))$.
\end{enumerate}
\end{theorem}

As already mentioned, the inclusion $span(T_\pi|\pi\in P(k,l))\subset Hom(u^{\otimes k},u^{\otimes l})$ is easy to check for $S_n$. The reverse inclusion can be justified by several methods, none of them being trivial. For instance this inclusion follows from the classical Tannaka duality.

In the free case now, we can use the embedding $S_n\subset S_n^+$, which at the level of Hom spaces gives reverse inclusions $Hom(u^{\otimes k},u^{\otimes l})\subset span(T_\pi|\pi\in P(k,l))$. Now since $S_n^+$ was obtained from $S_n$ by ``removing commutativity'', at the Tannakian level we should ``remove the basic crossing $X$'', and so we reach to the above sets $NC(k,l)$.

In this latter proof we have used of course the results in \cite{wo2}. See \cite{asn}.

We recall that the law of a self-adjoint variable $a\in (A,\varphi)$ is the real probability measure having as moments the numbers $\varphi(a^k)$. In what follows $\varphi$ will be the Haar functional.

\begin{theorem}
The law of the character $\chi=Tr(u)$ is as follows:
\begin{enumerate}
\item For $S_n$ with $n\to\infty$ we obtain the Poisson law, $p=\frac{1}{e}\sum_{r=0}^\infty\frac{\delta_r}{r!}$.

\item For $S_n^+$ with $n\geq 4$ we obtain the free Poisson law, $\pi=\frac{1}{2\pi}\sqrt{4x^{-1}-1}\,dx$.
\end{enumerate}
\end{theorem}

The first result is well-known, and has a very beautiful proof. Indeed, since the standard coordinates on $S_n$ are the characteristic functions $u_{ij}=\chi(\sigma|\sigma(j)=i)$, we obtain:
$$\chi=\sum_{i=1}^n\chi(\sigma|\sigma(i)=i)=|fix(.)|$$

Now by using the well-known fact that a random permutation has no fixed points with probability $1/e$, hence has $r$ fixed points with probability $1/(r!e)$, we get, as claimed:
$$\mathbb P(\chi(\sigma)=r)=\mathbb P(|fix(\sigma)|=r)\simeq\frac{1}{e}\cdot\frac{1}{r!}$$

Another proof, which works as well for $S_n^+$, is as follows. Let $D(k,l)$ be the sets of diagrams appearing in Theorem 2.3. By using basic Peter-Weyl theory, we have:
$$\int\chi^k=\dim(Fix(u^{\otimes k}))\simeq |D(0,k)|$$

Here the sign at right corresponds to our $n\to\infty$ and $n\geq 4$ assumptions, and comes from the fact that the maps $T_\pi$ are not exactly linearily independent. See \cite{asn}. Now, with this formula in hand, we conclude that the moments of $\chi$ are the Bell, respectively the Catalan numbers, and this gives the result.

Here is a purely functional analytic application of Theorem 2.4:

\begin{corollary}
For $n\geq 5$ the discrete quantum group $\widehat{S}_n^+$ is not amenable.
\end{corollary}

Indeed, recall the Kesten amenability criterion, which tells us that a discrete group $\Gamma=<g_1,\ldots,g_n>$, chosen with $\{g_i\}=\{g_i^{-1}\}$, is amenable if and only if:
$$n\in Spec(g_1+\ldots+g_n)$$

The point is that the Kesten operator is nothing but the trace of $u=diag(g_1,\ldots,g_n)$, which is the fundamental representation of $\widehat{\Gamma}$. Thus, the Kesten criterion reads:
$$n\in Spec(Tr(u))$$

With a bit more work, one can show that the same result holds for any discrete quantum group $\Gamma$. See \cite{sbf}. Thus, the above result follows from Theorem 2.4 (2) and from:
$$Supp(\pi)=[0,4]$$

There are many interesting questions arising from Corollary 2.5, and from the work of Brannan \cite{bra}, Fima \cite{fim}, Kyed \cite{kye} and Vaes-Vergnioux \cite{vve}. In all cases, the question is whether the techniques there apply to $S_n^+$, at $n\geq 5$.

Here is another corollary, this time of purely algebraic nature: 

\begin{corollary}
For $n\geq 4$ the irreducible representations of $S_n^+$ can be labeled by $\mathbb N$, and satisfy the Clebsch-Gordan rules: $r_a\otimes r_b=r_{|a-b|}+r_{|a-b|+1}+\ldots+r_{a+b}$.
\end{corollary}

This result can be indeed deduced from Theorem 2.4. The idea is that, with the initial data $r_0=1$ and $r_1=u-1$, we can construct a family of virtual representations $r_a$ satisfying the Clebsch-Gordan rules. The problem is to prove that these virtual representations $r_a$ are true representations, and are irreducible. But this can be done by recurrence on $a$, by using the fact that the moments of $\chi$ are the Catalan numbers:
$$\dim(Fix(u^{\otimes k}))=\frac{1}{k+1}\binom{2k}{k}$$

This proof is of course not very enlightening. Probably the best argument is as follows: let us enlarge attention to the groups $S(X)$ and quantum groups $S^+(X)$, where $X$ can be any finite noncommutative space. We have then the following result:

\begin{theorem}
If $u$ denotes the fundamental representation of $S^+(X)$ then
$$Hom(u^{\otimes k},u^{\otimes l})\simeq TL_n(k,l)$$
where $TL_n$ denotes the Temperley-Lieb algebra of index $n=|X|$.
\end{theorem}

As a first remark, for a usual space $X=\{1,\ldots,n\}$, this is exactly Theorem 2.3 (2) above, modulo a ``fattening'' of the diagrams. In general, the proof is based on the fact that the multiplication and unit of the algebra $C(X)$ can be represented as:
$$m=|\cup|\quad\quad u=\cap$$

Now since $S^+(X)$ is the universal quantum group acting on $C(X)$, at the Tannakian level we get that its associated category is simply given by: 
$$<m,u>=<|\cup|,\ \cap>=TL_n$$

With this conceptual result in hand, let us go back now to Corollary 2.6, and try to find a better proof. The point is that for the space $M_2=Spec(M_2(\mathbb C))$ we have:
$$S^+(M_2)=S(M_2)=SO_3$$

On the other hand, it is well-known that under the Jones index assumption $n\geq 4$, the fusion ring for the category $TL_n$ doesn't depend on $n$. Thus, the fusion rules for $S_n^+$ with $n\geq 4$ are the same as those for $SO_3$, as stated in Corollary 2.6 above.

\begin{problem}
Is it true that the compact quantum groups having Clebsch-Gordan fusion rules are exactly the quantum groups $S^+(X)$, with $|X|\geq 4$?
\end{problem}

The general method for attacking this type of problem is via Kazhdan-Wenzl reconstruction \cite{kwe}. In the type $B$ case the result here is due to Tuba and Wenzl \cite{twe}.

Observe that the above discussion shows that with $q<0$ given by $n=(q+q^{-1})^2$ we have a monoidal equivalence \cite{bdv}, i.e. an equivalence of tensor categories, as follows:
$$S_n^+\sim SO_3^q$$

Observe also that the dimensions of the representations of $S_n^+$ with $n\geq 5$ are ``exotic''. More precisely, if $x,y$ are the roots of $X^2-(n-2)X+1=0$, then:
$$dim(r_a)=\frac{x^{a+1}-y^{a+1}}{x-y}+\frac{x^a-y^a}{x-y}$$

There are as well a few other algebraic questions, coming from the work of Collins-H\"artel-Thom \cite{cht}, De Commer \cite{dec}, So\l tan \cite{sol}, Vergnioux \cite{vrg} and Voigt \cite{voi}.

Let us go back now to Theorem 2.4. One problem is to compute the value of the Haar functional on any element of type $u_{i_1j_1}\ldots u_{i_kj_k}$. The method here, developed for classical groups by Collins and \'Sniady in \cite{csn}, applies as well to $S_n,S_n^+$, and we have: 

\begin{theorem}
We have the Weingarten integration formula
$$\int u_{i_1j_1}\ldots u_{i_kj_k}=\sum_{\pi,\sigma\in D(k)}\delta_\pi(i)\delta_\sigma(j)W_{kn}(\pi,\sigma)$$
where $D(k)=D(0,k)$, and where $W_{kn}=G_{kn}^{-1}$, with $G_{kn}(\pi,\sigma)=n^{|\pi\vee\sigma|}$.
\end{theorem}

The proof of this formula, found in \cite{asn}, is based on Theorem 2.3. Consider indeed the $n^k\times n^k$ matrix $P$ obtained by integrating the coefficients of $u^{\otimes k}$. That is: 
$$P=\left(\int u_{i_1j_1}\ldots u_{i_kj_k}\right)_{i_1\ldots i_k,j_1\ldots j_k}$$

We know by Peter-Weyl theory that $P$ is the projection on $Fix(u^{\otimes k})$, and from Theorem 2.3 that $Fix(u^{\otimes k})$ is spanned by $D(k)$. The Gram matrix of $D(k)$ is:
$$<T_\pi,T_\sigma>=\sum_{i_1\ldots i_k}\delta_\pi(i)\delta_\sigma(i)=n^{|\pi\vee\sigma|}$$

By combining these ingredients and doing some linear algebra, this gives the result. As a first consequence, we have the following refinement of Theorem 2.4:

\begin{theorem}
The asymptotic laws of truncated characters $\chi_t=\sum_{i=1}^{[tn]}$ are:
\begin{enumerate}
\item For $S_n$: the Poisson law, $p_t=e^{-t}\sum_{r=0}^\infty\frac{t^r}{r!}\delta_r$.

\item For $S_n^+$: the free Poisson law, $\pi_t=\max (1-t,0)\delta_0+\frac{\sqrt{4t-(x-1-t)^2}}{2\pi x}\,dx$.
\end{enumerate}
\end{theorem}

This result is important for two reasons. First, it requires $n\to\infty$ for both $S_n,S_n^+$, hence fixes the $n\to\infty$ vs. $n\geq 4$ issue in Theorem 2.4. And second, it proves that the ``liberation'' operation $S_n\to S_n^+$ is compatible with the Bercovici-Pata bijection \cite{bpa}.

The proof uses the Weingarten formula. With $s=[tn]$, we get:
$$\int \chi_t^k=\sum_{\pi,\sigma\in D(k)}s^{|\pi\vee\sigma|}W_{kn}(\pi,\sigma)=Tr(G_{ks}W_{kn})$$

Now since with $n\to\infty$, the Gram matrix becomes equal to its diagonal, $G_{kn}=\Delta_{kn}+O(n^{-1})$, by inverting we get as well $W_{kn}=\Delta_{kn}^{-1}+O(n^{-1})$, and this gives:
$$\int \chi_t^k=Tr(\Delta_{ks}\Delta_{kn}^{-1})+O(n^{-1})=\sum_{\pi\in D(k)}t^{|\pi|}+O(n^{-1})$$

According now to the general theory in \cite{nsp}, this shows that the cumulants for $S_n$, as well as the free cumulants for $S_n^+$, are simply $t,t,t,\ldots$, and this gives the result.

Let $P_r$ be the Chebycheff polynomials, $P_0=1,P_1=n$ and $P_{r+1}=nP_r-P_{r-1}$, and:
$$a_k=\frac{1}{k+1}\binom{2k}{k},\quad\quad f_{kr}=\binom{2k}{k-r}-\binom{2k}{k-r-1}$$

With these notations, the we have the following results:

\begin{theorem}
The determinant of $G_{kn}$ is as follows:
\begin{enumerate}
\item For $S_n$: $\det(G_{kn})=\prod_{\pi\in P(k)}\frac{n!}{(n-|\pi|)!}$.

\item For $S_n^+$: $\det(G_{kn})=(\sqrt{n})^{a_k}\prod_{r=1}^kP_r(\sqrt{n})^{d_{kr}}$, with $d_{kr}=f_{kr}-f_{k+1,r}$.
\end{enumerate}
\end{theorem}

The classical result follows by an upper triangularization procedure. The free result uses Di Francesco's ``meander determinant'' in \cite{dgg}. See \cite{gra}.

We have as well the following generalization of Theorem 2.4:

\begin{theorem}
The asymptotic law of $\chi^{(r)}=Tr(u^r)$ is as follows:
\begin{enumerate}
\item For $S_n$: we obtain a linear combination of independent Poisson laws.

\item For $S_n^+$: $\chi^{(1)}$ is free Poisson, $\chi^{(2)}$ is semicircular, $\chi^{(r)}$ with $r\geq 3$ is circular.
\end{enumerate}
\end{theorem}

The classical result is due to Diaconis and Shahshahani \cite{dsh}. Observe that in terms of the eigenvalue list $\lambda_1\leq\ldots\leq\lambda_n$, the variable in the statement is:
$$Tr(u^r)=\lambda_1^r+\ldots+\lambda_n^r$$

In the free case the methods in \cite{dsh} cannot work. In addition, a surprise appears at $r=3$, where the variable in question is not self-adjoint:
$$Tr(u^3)=\sum_{ijk}u_{ij}u_{jk}u_{ki}\neq \sum_{ijk}u_{ki}u_{jk}u_{ij}=Tr(u^3)^*$$

We refer to our paper \cite{sas} for the proof for $S_n^+$, which applies as well to $S_n$.

\begin{question}
Is there a notion of ``eigenvalues'' for the quantum group $S_n^+$, which is compatible with the above results?
\end{question}

In other words, we would like to have a free analogue of the $n$ eigenvalue functions $\lambda_1\leq\ldots\leq\lambda_n$, making the formulae in Theorem 2.12 work. As observed in \cite{sas}, such eigenvalue functions can be constructed for certain related quantum groups $O_n^*,H_n^*$.

The main conceptual result regarding $S_n,S_n^+$ concerns their action on random variables. An action of $S=(S_n)$ or $S^+=(S_n^+)$ on a sequence of random variables $x_1,x_2,\ldots$ is a coaction on the generated algebra, leaving invariant the joint distribution.

\begin{theorem}
A sequence of random variables $x_1,x_2,\ldots$ is invariant iff:
\begin{enumerate}
\item $S$ case: $x_i$ are independent and i.d. with respect to the tail algebra.

\item $S^+$ case: $x_i$ are free and i.d. with respect to the tail algebra.
\end{enumerate}
\end{theorem}

Here the first result is the classical de Finetti theorem, and the second result is the free De Finetti theorem, due to K\"ostler-Speicher \cite{ksp}. A unified proof was given in \cite{def}.

\section{Actions on finite spaces}

We have seen in the previous section that the quantum groups $S_n^+$ are understood quite well. In this section we extend some of these results, to subgroups $G\subset S_n^+$.

Consider a finite graph $X$, with edges colored by real numbers. The basic examples are the finite metric spaces, and the usual (monocolored) graphs. If we label the vertices $1,\ldots,n$, the adjacency matrix becomes a real symmetric matrix, $d\in M_n(\mathbb R)$.

\begin{definition}
$G^+(X)$ be the biggest quantum group acting on $X$, in the sense that $C(G^+(X))$ is the quotient of $C(S_n^+)$ by the relations coming from $du=ud$.
\end{definition}

As a first remark, the classical version of $G^+(X)$ is the usual symmetry group $G(X)$. Indeed, with $u_{ij}=\chi(\sigma|\sigma(j)=i)$, the relation $du=ud$ reads:
$$d_{ij}=d_{\sigma(i)\sigma(j)}$$

Let us look now more carefully at the relation $du=ud$. Since $d$ is a real symmetric matrix, we can consider its spectral decomposition:
$$d=\sum_\lambda\lambda\cdot P_\lambda$$

Also, for any $r\in\mathbb R$ we can consider the matrix $d^r$ given by $d^r_{ij}=1$ if $d_{ij}=r$, and $d^r_{ij}=0$ if not. We obtain in this way the ``color decomposition'' of $d$:
$$d=\sum_rr\cdot d^r$$

Here, as in the case of the spectral decomposition, the sum is over nonzero terms.

\begin{proposition}
For a quantum group $G\subset S_n^+$, the following are equivalent:
\begin{enumerate}
\item We have $du=ud$, i.e. $G\subset G^+(X)$.

\item We have $P_\lambda u=uP_\lambda$, for any $\lambda\in\mathbb R$.

\item We have $d^ru=ud^r$, for any $r\in\mathbb R$.
\end{enumerate}
\end{proposition}

Here the equivalence between (1) and (2) follows from definitions, and the equivalence with (3) follows by doing a bit of matrix analysis. These results are very useful, because they can be further combined: for instance $u$ must commute with any color component $P_\lambda^r$ of any spectral projection $P_\lambda$, and so on. This method of ``spectral-color'' decomposition proves to be quite efficient for the computation of $G^+(X)$. See \cite{hgr}.

Observe that the equivalence between (1) and (3) can be reformulated as follows, where $X^r$ denotes the usual (monocolored) graph having adjacency matrix $d^r$:
$$G^+(X)=\bigcap_rG^+(X^r)$$

More generally now, let $X$ be a finite noncommutative set, endowed with a self-adjoint operator $d:L^2(X)\to L^2(X)$, where the $L^2$ space is with respect the counting measure.

\begin{definition}
$G^+(X)$ be the biggest quantum group acting on $X$, in the sense that $C(G^+(X))$ is the quotient of $C(S^+(X_{set}))$ by the relations coming from $du=ud$.
\end{definition}

Here we have used the quantum groups $S^+$ appearing in Theorem 2.7 above. This definition is of course very general, and covers for instance all the 0-dimensional spectral triples in the sense of Connes \cite{con}. We will be back to these examples, a bit later.

Let $*$ be the free product of discrete quantum groups \cite{wa1}:
$$C^*(\Gamma*\Lambda)=C^*(\Gamma)*C^*(\Lambda)$$

We denote by $\hat{*}$ the dual operation, at the level of compact quantum groups:
$$C(G\;\hat{*}\;H)=C(G)*C(H)$$

Finally, let us denote by $X_1\sqcup\ldots\sqcup X_k$ a disconnected union of graphs.

\begin{proposition}
We have $G^+(X_1)\;\hat{*}\;\ldots\;\hat{*}\;G^+(X_k)\subset G^+(X_1\sqcup\ldots\sqcup X_k)$.
\end{proposition}

This follows from definitions, and might seem to be yet another triviality. However, let us work out the simplest example. For the graph formed by two segments we get: 
$$\mathbb Z_2\;\hat{*}\;\mathbb Z_2\subset G^+(|\ \ |)$$

The dual free product on the left is by definition the dual of $\mathbb Z_2*\mathbb Z_2=D_\infty$. As for the quantum group of the right, we prefer to view it as the quantum symmetry group of the complement of $X=|\ \ |$, which is $X^c=\square$. Thus, the above inclusion reads:
$$\widehat{D}_\infty\subset G^+(\square)$$

And this reminds of course the inclusion $\widehat{D}_\infty\subset S_4^+$ from Proposition 1.7 above, which was a bit the start of everything! More precisely, we have reached now to the following quite conceptual conclusion: ``$S_4^+$ is infinite simply because $\widehat{D}_\infty$ acts on the square''.

Consider now an arbitrary union $kX=X\sqcup\ldots\sqcup X$ ($k$ terms).

\begin{proposition}
If $X$ is connected then $G^+(kX)=G^+(X)\wr_*S_k^+$.
\end{proposition}

Here we use the free wreath product operation, constructed by Bichon in \cite{bi1}:
$$C(G\wr_*H)=(C(G)^{*k}*C(H))/<[u_{ij}^{(a)},v_{ab}]=0>$$

More generally now, let $X,Y$ be finite graphs, and consider the graphs $X*Y$, $X\circ Y$ and $X\times Y$ having as vertex set $X_{set}\times Y_{set}$, and with the following edges:
$$(i,\alpha)\sim(j,\beta)\iff
\begin{cases}
\alpha\sim\beta\ {\rm or}\ \alpha=\beta,\,i\sim j&(*\ {\rm case})\\
i\sim j,\,\alpha\sim\beta&(\circ\ {\rm case})\\
i=j,\,\alpha\sim\beta\ {\rm or}\ i\sim j,\alpha=\beta&(\times\ {\rm case})
\end{cases}$$
 
In terms of the adjacency matrix the formulae are as follows, where $\mathbb I$ denotes the square matrix filled with $1$'s:
$$d=
\begin{cases}
d_X\otimes 1+\mathbb I\otimes d_Y&(*\ {\rm case})\\
d_X\otimes d_Y&(\circ\ {\rm case})\\
d_X\otimes 1+1\otimes d_Y&(\times\ {\rm case})
\end{cases}$$

Let us first look at $*$. This product is obtained by putting a copy of $X$ at each vertex of $Y$, so in particular we have $kX=X*Y_k$, where $Y_k$ is the $k$-point graph.

\begin{theorem}
Let $X,Y$ be regular graphs, with spectra $\{\lambda_i\},\{\mu_k\}$. If $X$ is connected and
$$\{\lambda_1-\lambda_i|i\neq 1\}\cap\{-n\mu_k\}=\emptyset$$
where $n$ and $\lambda_1$ are the order and valence of $X$, then $G^+(X*Y)=G^+(X)\wr_*G^+(Y)$.   
\end{theorem}

This result, generalizing Proposition 3.5, was proved in \cite{qgr}. Now let us look at $\circ,\times$. Under suitable assumptions on the eigenvalues of $X,Y$, we will get the spectral decomposition of $d$ in terms of that of $d_X,d_Y$. If we assume in addition that $X,Y$ are regular and connected, the 1-dimensional eigenspaces coming from the valence can be used for ``splitting'' the action of $G^+$, and we have the following result, once again from \cite{qgr}:

\begin{theorem}
Let $X,Y$ be connected regular graphs, with spectra $\{\lambda_i\},\{\mu_k\}$.
\begin{enumerate}
\item If $\{\lambda_i/\lambda_j\}\cap\{\mu_k/\mu_l\}=\{1\}$ then $G^+(X\circ Y)=G^+(X) \times G^+(Y)$. 

\item If $\{\lambda_i-\lambda_j\}\cap\{\mu_k-\mu_l\}=\{0\}$ then $G^+(X\times Y)=G^+(X)\times G^+(Y)$.
\end{enumerate}
\end{theorem}

With these results in hand, it is tempting to try to compute $G^+(X)$ for all vertex-transitive graphs of small order. Here is the result found in \cite{qgr}, up to order 9:

\begin{theorem}
The vertex-transitive graphs of order $\leq 9$, modulo complementation, and their quantum symmetry groups, are as follows:
\begin{enumerate}
\item Unions of simplices $rK_n$, having $G^+=S_n^+\wr_*S_r^+$.

\item Cycles $C_n$, $n\neq 4$, and the cycles with chords $C_8^2,C_9^3$, having $G^+=D_n$.

\item The cube, having $G^+=\mathbb Z_2\times S_4^+$, and $K_3\times K_3$, having $G^+=S_3\wr\mathbb Z_2$.

\item The two squares, $2C_4$, having $G^+=(\mathbb Z_2\wr_*\mathbb Z_2)\wr_*\mathbb Z_2$.
\end{enumerate}
\end{theorem}

We actually have a longer result in \cite{qgr}, going up to order 11, but not covering the Petersen graph. This well-known graph, having 10 vertices, consists of a pentagon with a pentagram inside, with 5 spokes. Our main question, subsequent to \cite{qgr}, is of course:

\begin{problem}
What is the quantum symmetry group of the Petersen graph?
\end{problem}

Regarding now the proof of Theorem 3.8, this basically follows from Theorem 3.6 and Theorem 3.7, and from a number of ``no quantum symmetry'' results, coming from Theorem 3.6, from Theorem 3.7, or (e.g. for $K_3\times K_3$) from a long algebraic computation based on Proposition 3.2. The subject is quite mysterious, for instance because:
$$G(C_n)=G^+(C_n)\iff n\neq 4$$

Indeed, the $\implies$ assertion follows from the discussion after Proposition 3.4. As for the converse, this follows from Proposition 3.2, after a brief computation. See \cite{hgr}.

There is an obvious link here with the following conjecture of Goswami \cite{go2}:

\begin{conjecture}
A compact connected Riemannian manifold cannot have genuine quantum isometries.
\end{conjecture}

More precisely, Goswami associated in \cite{go1} a quantum isometry group $G^+(M)$ to any noncommutative compact Riemannian manifold $M$ in the sense of Connes \cite{con}. In the 0-dimensional case the idea is basically the one in Definition 3.3 above, and in the general case the idea is similar, but more technical, with $d$ being the Laplacian of $M$. 

Goswami was mostly interested in genuine noncommutative manifolds, like the Standard Model one, whose quantum isometries have several potential applications \cite{bdd}. However, his construction makes sense for any manifold, and leads to the above conjecture.

Back to graphs now, a very interesting example is the $n$-cube. The symmetry theory of this graph, heavily investigated in \cite{ahn}, can be summarized as follows:

\begin{theorem}
Let $X=[-1,1]^n$ be the cube in $\mathbb R^n$, regarded as graph, and let $Y$ be the intersection of $X$ with the $n$ coordinate axes, regarded as a graph formed by $n$ segments.
\begin{enumerate}
\item $G(X)=G(Y)=\mathbb Z_2\wr S_n$.

\item $G^+(X)=O_n^{-1}$.

\item $G^+(Y)=\mathbb Z_2\wr_*S_n^+$
\end{enumerate}
\end{theorem}

In this statement the group in (1) is the hyperoctahedral group $H_n$. The quantum group in (2) is a Drinfeld-Jimbo twist at $q=-1$. This twist has the same fusion rules as $O_n$, so the operation $H_n\to O_n^{-1}$ is not a ``liberation'', in the sense of Theorem 2.10.

As for the quantum group in (3), which obviously looks as the good generalization of $H_n$, this is called the hyperoctahedral quantum group, and is denoted $H_n^+$.

This fact that $H_n\to H_n^+$ is indeed a ``liberation'' in the probabilistic sense comes from the following result, once again from \cite{ahn}:

\begin{theorem}
The asymptotic laws of truncated characters $\chi_t=\sum_{i=1}^{[tn]}u_{ii}$ are:
\begin{enumerate}
\item For $H_n$: the Bessel law, $b_t=law(\alpha-\beta)$, where $\alpha,\beta\sim p_{t/2}$ are independent.

\item For $H_n^+$: the free Bessel law, $\beta_t=law(\alpha-\beta)$, where $\alpha,\beta\sim\pi_{t/2}$ are free.
\end{enumerate}
\end{theorem}

Observe the similarity with Theorem 2.10. The terminology comes from the fact that the density of $b_t$ is given by a Bessel function of the first kind:
$$b_t=e^{-t}\sum_{k=-\infty}^\infty\left(\sum_{p=0}^\infty \frac{(t/2)^{|k|+2p}}{(|k|+p)!p!}\right)\delta_k$$

As for the measure $\beta_t$, this has remarkable properties as well. For instance the odd moments vanish, and the even moments involve the Fuss-Narayana numbers:
$$\int x^{2k}\,d\beta_t(x)=\sum_{b=1}^k\frac{1}{b}\begin{pmatrix}k-1\cr b-1\end{pmatrix}\begin{pmatrix}2k\cr b-1\end{pmatrix}t^b$$

At the level of diagrams, the partitions as in Theorem 2.3 are those having even blocks, and the algebra as in Theorem 2.7 is the Fuss-Catalan algebra of Bisch and Jones \cite{bj1}.

These results were subsequently explored and generalized to the case of the complex reflection groups $H_n^s=\mathbb Z_s\wr S_n$ and their free versions $H_n^{s+}=\mathbb Z_s\wr_*S_n^+$. See \cite{fbl}, \cite{fru}.

A general diagrammatic formalism, covering the quantum groups $S_n,S_n^+,H_n,H_n^+$, as well as many other examples, such as the orthogonal group $O_n$ and its free analogue $O_n^+$ constructed by Wang in \cite{wa1}, was developed in \cite{olg}, then in \cite{cls}, \cite{sas}, \cite{def}:

\begin{definition}
A compact quantum group $S_n\subset  G\subset O_n^+$ is called ``easy'' if
$$Hom(u^{\otimes k},u^{\otimes l})=span(T_\pi|\pi\in D(k,l))$$
for certain sets of partitions $D(k,l)\subset P(k,l)$.
\end{definition}

This special class of quantum groups, very useful for unifying certain combinatorial or probabilistic results, was further studied in \cite{cur}, \cite{csp}, and in \cite{rau}.

A promising direction here comes from noncommutative geometry. First, the quantum isometry groups of 0-dimensional manifolds, generalizing the graph construction, were investigated in \cite{bgs}. In the higher dimensional case, a lot of interesting examples come from discrete group duals \cite{bsk}. The recent results in \cite{sym} show that some of these latter quantum isometry groups are ``super-easy'', in the sense that they appear as above, but with the representation $\pi\to T_\pi$ altered by an involution $J:\mathbb C^n\to\mathbb C^n$. 

Back now to $S_n^+$, the free hypergeometric law of parameters $(n,m,N)$ is the law of:
$$X(n,m,N)=\sum_{i=1}^n\sum_{j=1}^mu_{ij} \in C(S_N^+)$$

The terminology comes from the fact that the variable $X'(n,m,N)$, defined as above, but over the algebra $C(S_N)$, follows a hypergeometric law of parameters $(n,m,N)$.

\begin{theorem}
The moments of $X(n,n,n^2)$ are given by
$$\int X(n,n,n^2)^k\,dx=\frac{n^k}{(n+1)^k}\cdot\frac{q+1}{q-1}\cdot\frac{1}{k+1}\sum_{r=-k-1}^{k+1}(-1)^r\begin{pmatrix}2k+2\cr k+r+1\end{pmatrix}\frac{r}{1+q^r}$$
where $q\in [-1,0)$ is given by $q+q^{-1}=-n$.
\end{theorem}

The proof of this result, from \cite{tga}, and heavily based on \cite{sfo}, first relies on a purely algebraic twisting result, followed by a monoidal equivalence:
$$(S_{n^2}^+)^\tau=PO_n^+\sim SO^q_3$$

The remaining computation can be done by using a number of delicate manipulations involving Askey-Wilson polynomials \cite{awi}, and leads to the above formula.

These computations are of great interest, for instance in connection with the Atiyah conjecture for quantum groups \cite{cht}, and with the free spheres introduced in \cite{ncs}.

We would like to end this section with a key algebraic question: 

\begin{problem}
Is the inclusion $S_n\subset S_n^+$ maximal?
\end{problem}

The answer here is yes at $n=1,2,3$, and also yes at $n=4$, due to the results in \cite{qfp}. It is also yes in the ``easy'' case, as shown in \cite{cls}. There are some similarities here with a recent result in \cite{max}, stating that a certain related inclusion $O_n\subset O_n^*$ is maximal.

\section{Hadamard matrices}

A complex Hadamard matrix is a matrix $H\in M_n(\mathbb C)$ whose entries are on the unit circle, $|H_{ij}|=1$, and whose rows are pairwise orthogonal.

It follows from definitions that the columns of $H$ are pairwise orthogonal as well. 

The link with the quantum permutation groups is immediate. Let $H_1,\ldots,H_n$ be the rows of $H$, regarded as invertible elements of the algebra $\mathbb C^n$.

\begin{proposition}
For an Hadamard matrix $H\in M_n(\mathbb C)$, the rank one projections
$$P_{ij}=Proj(H_i/H_j)$$
form a magic unitary matrix $P\in M_n(A)$, where $A=M_n(\mathbb C)$.
\end{proposition}

Indeed, the orthogonality condition on the rows of $P$ follows from the fact that the corresponding vectors form an orthogonal basis of $\mathbb C^n$:
$$\Big\langle\frac{H_i}{H_j},\frac{H_i}{H_k}\Big\rangle=\sum_r\bar{H}_{jr}H_{kr}=\delta_{jk}$$

As for the orthogonality condition on the columns, this follows from a similar computation. According now to Definition 1.4 above, we have:

\begin{proposition}
Associated to any Hadamard matrix $H\in M_n(\mathbb C)$ is the morphism
$$\pi_H:C(S_n^+)\to M_n(\mathbb C)$$
given by $u_{ij}\to Proj(H_i/H_j)$, where $H_1,\ldots,H_n$ are the rows of $H$.
\end{proposition}

We will be back to this morphism in section 5 below. The idea there will be that associated to $\pi_H$ is a quantum group $G\subset S_n^+$, describing the ``symmetries'' of $H$.

For the moment, let us try to understand what the examples are. First, we have the Fourier matrix, based on the root of unity $w=e^{2\pi i/n}$:
$$F_n=\begin{pmatrix}
1&1&1&\ldots&1\\
1&w&w^2&\ldots&w^{n-1}\\
\ldots&\ldots&\ldots&\ldots&\ldots\\
1&w^{n-1}&w^{2(n-1)}&\ldots&w^{(n-1)^2}
\end{pmatrix}$$

This matrix is the matrix of the Fourier transform over $\mathbb Z_n$. There are many other examples of Hadamard matrices, usually coming from finite groups, or other discrete structures. See \cite{tz2}. An interesting observation, due to Popa, is that a unitary $U\in U_n$ satisfies $\Delta\perp U\Delta U^*$ if and only if $H=\sqrt{n}\cdot U$ is Hadamard. See \cite{pop}.

In order to present now some classification results, we use the following definitions:

\begin{definition}
Let $H,K$ be two complex Hadamard matrices.
\begin{enumerate}
\item $H$ is called dephased if its first row and column consist only of $1$'s.

\item $H,K$ are called equivalent if one can pass from one to the other by permuting the rows or columns, or by multiplying them by complex numbers of modulus $1$.
\end{enumerate}
\end{definition}

Observe that any complex Hadamard matrix can be supposed to be in dephased form, up to the above equivalence relation. With a few exceptions, we will always do so.

At $n=2,3$ the Fourier matrix is the only one, up to equivalence. At $n=4$ we have the following family, depending on a parameter on the unit circle, $q\in\mathbb T$:
$$F_4^q=\begin{pmatrix}
1&1&1&1\\
1&-1&q&-q\\
1&1&-1&-1\\ 
1&-1&-q&q
\end{pmatrix}$$

Observe that at $q=1$ we have $F_4^q=F_2\otimes F_2$. Observe also that at $q=-1$ we have $F_4^q\simeq F_2\otimes F_2$, and that at $q=\pm i$ we have $F_4^q\simeq F_4$. We will come back a bit later to these observations, with a discussion regarding the deformations of $F_n$.

The following remarkable result is due to Haagerup \cite{haa}.

\begin{theorem}
The only Hadamard matrices at $n\leq 5$ are $F_2,F_3,F_4^q,F_5$.
\end{theorem}

At $n=6$ we have several deformations of $F_6=F_2\otimes F_3=F_3\otimes F_2$. These matrices appear as particular cases of the following general construction, due to Di\c t\u a \cite{dit}:

\begin{definition}
The Di\c t\u a deformation of a tensor product $H\otimes K\in M_{nm}(\mathbb C)$, with matrix of parameters $L\in M_{m\times n}(\mathbb T)$, is $H\otimes_LK=(H_{ij}L_{aj}K_{ab})_{ia,jb}$.
\end{definition}

As a first observation, by using the above equivalence relation, we can always assume that $L$ is ``dephased'', in the sense that its first row and column consist of $1$'s.

The basic example comes from deforming the tensor product $F_2\otimes F_2$, with the matrix of parameters $L=(^1_1{\ }^1_q)$. We obtain in this way the above matrix $F_4^q$:
$$F_4^q=
\begin{pmatrix}
1&1\\
1&-1
\end{pmatrix}
\otimes_{\begin{pmatrix}
1&1\\
1&q
\end{pmatrix}}
\begin{pmatrix}
1&1\\
1&-1
\end{pmatrix}$$

As already mentioned, this matrix can be thought of as being a ``deformation'' of $F_4$, because at $q=\pm i$ we have $F_4^q\simeq F_4$. More generally, we have:

\begin{proposition}
$F_{nm}=F_n\otimes_LF_m$, where $L=(w^{(a-1)(j-1)})_{aj}$ with $w=e^{2\pi i/nm}$.
\end{proposition}

Let us go back now to the enumeration of Hadamard matrices of small order. At $n=6$ we have the (dephased) Di\c t\u a deformations of $F_2\otimes F_3$ and $F_3\otimes F_2$, given by:
$$F_6^{(^r_s)}=F_2
\otimes_{\begin{pmatrix}
1&1\\
1&r\\
1&s
\end{pmatrix}}
F_3\,,\quad\quad\quad
F_6^{(rs)}=F_3
\otimes_{\begin{pmatrix}
1&1&1\\
1&r&s
\end{pmatrix}}
F_2$$

Here $r,s$ are two parameters on the unit circle, $r,s\in\mathbb T$. In matrix form:
$$F_6^{(^r_s)}=\begin{pmatrix}
1&1&1&1&1&1\\
1&j&j^2&r&jr&j^2r\\
1&j^2&j&s&j^2s&js\\ 
1&1&1&-1&-1&-1\\
1&j&j^2&-r&-jr&-j^2r\\
1&j^2&j&-s&-j^2s&-js
\end{pmatrix},\quad
F_6^{(rs)}
=\begin{pmatrix}
1&1&1&1&1&1\\
1&-1&r&-r&s&-s\\
1&1&j&j&j^2&j^2\\ 
1&-1&jr&-jr&j^2s&-j^2s\\
1&1&j^2&j^2&j&j\\
1&-1&j^2r&-j^2r&js&-js
\end{pmatrix}$$

There are many other examples at $n=6$, and we will restrict attention to a certain special class of matrices, defined as follows:

\begin{definition}
A complex Hadamard matrix is called regular if all the scalar products between distinct rows decompose as sums of cycles.
\end{definition}

Here by ``cycle'' we mean a sum of $p$-roots of unity, rotated by a scalar $\lambda\in\mathbb T$:
$$C=\sum_{k=1}^p\lambda e^{2k\pi i/p}$$

Besides the above matrices, we have the Haagerup and Tao matrices \cite{haa}, \cite{tao}:
$$H^q=\begin{pmatrix}
1&1&1&1&1&1\\
1&-1&i&i&-i&-i\\ 
1&i&-1&-i&q&-q\\ 
1&i&-i&-1&-q&q\\
1&-i&\bar{q}&-\bar{q}&i&-1\\ 
1&-i&-\bar{q}&\bar{q}&-1&i
\end{pmatrix},\quad 
T=\begin{pmatrix}
1&1&1&1&1&1\\ 
1&1&j&j&j^2&j^2\\ 
1&j&1&j^2&j^2&j\\
1&j&j^2&1&j&j^2\\ 
1&j^2&j^2&j&1&j\\ 
1&j^2&j&j^2&j&1
\end{pmatrix}$$

The following result was proved in \cite{roq}:

\begin{theorem}
The only regular matrices at $n=6$ are $F_6^{(^r_s)},F_6^{(rs)},H^q,T$.
\end{theorem}

In the non-regular case, the general classification problem is open at $n=6$. See \cite{tz1}. 

One interesting example at $n=6$, with circulant structure, is the Bj\"orck-Fr\"oberg matrix \cite{bfr}, built by using one of the two roots of $a^2-(1-\sqrt{3})a+1=0$:
$$BF=\begin{pmatrix}
1&ia&-a&-i&-\bar{a}&i\bar{a}\\
i\bar{a}&1&ia&-a&-i&-\bar{a}\\
-\bar{a}&i\bar{a}&1&ia&-a&-i\\
-i&-\bar{a}&i\bar{a}&1&ia&-a\\
-a&-i&-\bar{a}&i\bar{a}&1&ia\\
ia&-a&-i&-\bar{a}&i\bar{a}&1
\end{pmatrix}$$

In a remarkable paper \cite{bni}, Beauchamp and Nicoara classified all self-adjoint Hadamard matrices at $n=6$. Their proof, and the proof of Theorem 4.8 above too, is heavily combinatorial. At $n=7$ we have the following matrix, discovered by Petrescu \cite{pet}:
$$P^q
=\begin{pmatrix}
1&1&1&1&1&1&1\\
1&qw&qw^4&w^5&w^3&w^3&w\\
1&qw^4&qw&w^3&w^5&w^3&w\\
1&w^5&w^3&\bar{q}w&\bar{q}w^4&w&w^3\\
1&w^3&w^5&\bar{q}w^4&\bar{q}w&w&w^3\\
1&w^3&w^3&w&w&w^4&w^5\\
1&w&w&w^3&w^3&w^5&w^4
\end{pmatrix}$$

Here $w=e^{2\pi i/6}$. This matrix, a non-trivial deformation of prime order, was found by using a computer program, and came as a big surprise at the time of \cite{pet}.

We can see from the above examples that, at least in the regular case, the roots of unity play a key role in the construction of complex Hadamard matrices.

\begin{definition}
The level of a complex Hadamard matrix $H\in M_n(\mathbb C)$ is the smallest number $l\in\{2,3,\ldots,\infty\}$ such that all the entries of $H$ are $l$-roots of unity.
\end{definition}

Here we agree that a root of unity of order $l=\infty$ is simply an element on the unit circle. It is convenient to introduce as well the following notions:

\begin{definition}
The Butson class $H_n(l)$ consists of Hadamard matrices in $M_n(\mathbb C)$ having as entries the $l$-th roots of unity. In particular:
\begin{enumerate}
\item $H_n(2)$ is the set of all $n\times n$ real Hadamard matrices.

\item $H_n(l)$ is the set of $n\times n$ Hadamard matrices of level $l'|l$.

\item $H_n(\infty)$ is the set of all $n\times n$ complex Hadamard matrices.
\end{enumerate}
\end{definition}

As a first example, let $H\in H_n(2)$, with $n\geq 3$. By using the equivalence relation, we may assume that the first three rows look blockwise as follows:  
$$H=\begin{pmatrix}
1&1&1&1\\
1&1&-1&-1\\
1&-1&1&-1\\
\ldots&\ldots&\ldots&\ldots
\end{pmatrix}$$

Now let $a,b,c,d$ be the lengths of the blocks in the third row. The orthogonality relations between the first three rows give $a+b=c+d$, $a+c=b+d$ and $a+d=b+c$, so we have $a=b=c=d$, and we can conclude that we have $4|n$.

The Hadamard conjecture states that the converse is true: 

\begin{conjecture}[Hadamard conjecture]
If $4|n$ then $H_n(2)\neq\emptyset$.
\end{conjecture}

This question, going back to the 19th century, is reputed to be of remarkable difficulty. The numeric verification so far goes up to $n=664$. See \cite{kta}.

This conjecture has as well an analytic interpretation, that we would like now to explain. Pick any $U\in U_n$ and apply the Cauchy-Schwarz inequality to its entries:
$$\sum_{i,j=1}^n1\cdot |U_{ij}|\leq\left(\sum_{i,j=1}^n1^2\right)^{1/2}\left(\sum_{i,j=1}^n|U_{ij}|^2\right)^{1/2}$$

We recognize at left the 1-norm of $U$. The right term being $n\sqrt{n}$, we get:
$$||U||_1\leq n\sqrt{n}$$

The point now is that the equality holds when the numbers $|U_{ij}|$ are proportional, and since the sum of squares of these numbers is $n^2$, we conclude that we have equality if and only if $|U_{ij}|=1/\sqrt{n}$, i.e. if and only if $H=\sqrt{n}\cdot U$ is Hadamard. Thus we have:

\begin{proposition}
For $U\in U_n$ we have $||U||_1\leq n\sqrt{n}$, with equality if and only if $U=H/\sqrt{n}$, for a certain complex Hadamard matrix $H$.
\end{proposition}

Now since the maximum of any positive function, and in particular of $U\to||U||_1$, can be obtained from moments, this suggests considering the following type of integral: 
$$I_G=\left(\int_G||U||_1^k\,dU\right)^{1/k}$$ 

More precisely, the Hadamard conjecture can be restated as follows:
$$4|n\implies\lim_{k\to\infty}I_{O_n}=n\sqrt{n}$$

This observation was made in \cite{omx}, and a number of related results were obtained there. For instance, it was shown that the Hadamard conjecture is equivalent to:
$$4|n\implies\lim_{k\to\infty}I_{O_n}\geq n\sqrt{n}-\frac{1}{n\sqrt{n}}$$

This estimate is almost optimal, so the Hadamard conjecture asks for the computation of $I_{O_n}$ with very high accuracy. In the general context of $O_n$ integration, where most results are in the $n\to\infty$ limit, this kind of problematics is quite new. See \cite{int}. 

We have the following theoretical questions, regarding the above integrals:

\begin{problem}
Let $I_G=(\int_G||U||_1^k\,dU)^{1/k}$.
\begin{enumerate}
\item Does the first-order term of $I_{O_n}$ count $|H_n(2)|$?

\item Is it possible to generalize this to any $l<\infty$?

\item What does the first-order term of $I_{U_n}$ count?
\end{enumerate}
\end{problem}

These questions are of course very vague, because we don't know what ``first-order term'' should mean. As an answer to the third question, we would expect of course some very simple numeric invariant of the complex Hadamard matrix manifold $H_n(\infty)$.

Let us go back now to the arbitrary sets of Butson matrices, $H_n(l)$ with $l<\infty$.

\begin{theorem}
We have the following list of obstructions:
\begin{enumerate}
\item Lam-Leung: If $H_n(l)\neq\emptyset$ and $l=p_1^{a_1}\ldots p_s^{a_s}$ then $n\in p_1\mathbb N+\ldots+p_s\mathbb N$.

\item de Launey: If $H_n(l)\neq\emptyset$ then there is $d\in\mathbb Z[e^{2\pi i/l}]$ such that $|d|^2=n^n$.

\item Sylvester: If $H_n(2)\neq\emptyset$ then $n=2$ or $4|n$.

\item Haagerup: If $H_5(l)\neq\emptyset$ then $5|l$.
\end{enumerate}
\end{theorem}

In this statement (1) follows from \cite{lle}, (2), where $d=\det H$, was studied in \cite{lau}, (3) was already explained, and (4) comes from Theorem 4.4 above.

The proof of the Sylvester obstruction, presented above, can be adapted for higher values of $l$, and leads to the following supplementary ``Sylvester obstructions'':
\begin{enumerate}
\item If $H_n(l)\neq\emptyset$ and $n=p+2$ with $p\geq 3$ prime, then $l\neq 2p^b$.

\item If $H_n(l)\neq\emptyset$ and $n=2q$ with $p>q\geq 3$ primes, then $l\neq 2^ap^b$.
\end{enumerate}

These results were proved in \cite{roq}. The point now is that, with these results in hand, the case $n\leq 10,l\leq 14$ is fully covered. We have the following table from \cite{roq}, describing for each $n,l$ either an explicit matrix in $H_n(l)$, or an obstruction which applies:
\setlength{\extrarowheight}{2pt}
{\small\begin{center}
\begin{tabular}[t]{|l|l|l|l|l|l|l|l|l|l|l|l|l|l|l|l|}
\hline $n\backslash l$\!
&2&3&4&5&6&7&8&9&10&11&12&13&14\\
\hline 2&$F_2$&$\circ$&$F_2$&$\circ$&$F_2$&$\circ$&$F_2$&$\circ$&$F_2$&$\circ$&$F_2$&$\circ$&$F_2$\\
\hline 3&$\circ$&$F_3$&$\circ$&$\circ$&$F_3$&$\circ$&$\circ$&$F_3$&$\circ$&$\circ$&$F_3$&$\circ$&$\circ$\\
\hline 4&$F_{22}$&$\circ$&$F_{22}$&$\circ$&$F_{22}$&$\circ$&$F_{22}$&$\circ$&$F_{22}$&$\circ$&$F_{22}$&$\circ$&$F_{22}$\\
\hline 5&$\circ$&$\circ$&$\circ$&$F_5$&$\circ_l$&$\circ$&$\circ$&$\circ$&$F_5$&$\circ$& $\circ_h$&$\circ$&$\circ$\\
\hline 6&$\circ_s$&$T$&$H^1$&$\circ$&$T$&$\circ$&$H^1$&$T$&$\circ_{s}$&$\circ$&$T$&$\circ$&$\circ_{s}$\\
\hline 7&$\circ$&$\circ$&$\circ$&$\circ$&$P^1$&$F_7$&$\circ$&$\circ$&$\circ_{s}$&$\circ$&$P^1$&$\circ$&$F_7$\\
\hline 8&$F_{222}\!$&$\circ$&$F_{222}\!$&$\circ$&$F_{222}\!$&$\circ$&$F_{222}\!$&$\circ$&$F_{222}\!$&$\circ$&$F_{222}\!$&$\circ$&$F_{222}\!$\\
\hline 9&$\circ$&$F_{33}$&$\circ$&$\circ$&$F_{33}$&$\circ$&$\circ$&$F_{33}$&$W$&$\circ$&$F_{33}$&$\circ$&$\circ_{s}$\\
\hline 10\!&$\circ_s$&$\circ$&$X$&$Y$&$Z$&$\circ$&$X$&$\circ$&$F_{10}$&$\circ$&$X$&$\circ$&$\circ$\\
\hline
\end{tabular}
\end{center}}
\setlength{\extrarowheight}{0pt}
\bigskip

Here $\circ,\circ_l,\circ_s,\circ_h$ denote the various obstructions in Theorem 4.14 and in the discussion afterwards, for $k_1,\ldots,k_s\in\{2,3\}$ we use the notation $F_{k_1\ldots k_s}=F_{k_1}\otimes\ldots\otimes F_{k_s}$, and $X,Y,Z,W$ are matrices constructed by a computer, and presented in \cite{roq}.

\begin{question}
What happens at $n=11$? And, what happens at $l=15$?
\end{question}

We refer here to \cite{roq} for technical details.

Finally, let us mention the following key conjecture, from \cite{roq}.

\begin{conjecture}
Any Butson matrix is regular.
\end{conjecture}

One difficulty in approaching this conjecture lies in the fact that the vanishing sums of roots of unity have a quite complicated combinatorics, cf. \cite{lle}.

\section{The correspondence}

We have seen in the previous section that the complex Hadamard matrices, while strongly related to the quantum permutation groups via tha magic condition, seem rather to belong to another planet. There is no obvious relation between the substantial mass of combinatorics presented in sections 1-3 above, and the substantial mass of combinatorics presented in section 4 above. So, the following statement might seem a bit strange:

\begin{fact}
The symmetries of a complex Hadamard matrix $H\in M_n(\mathbb C)$ are described by a quantum permutation group $G\subset S_n^+$. Conversely, each quantum permutation group $G\subset S_n^+$ is the symmetry group of an Hadamard-type object $H\in\tilde{M}_n(\mathbb C)$.
\end{fact}

In this section we will try to explain why this fact is true. The fact that the known combinatorics of quantum permutation groups remains quite far from the known combinatorics of complex Hadamard matrices will be discussed in the next section. We will propose there some conjectural answers to the question, making the above fact a potentially useful tool, in the study of both quantum permutations and Hadamard matrices.

Let us begin with the notion of Hopf image, recently axiomatized in \cite{ifr}:

\begin{definition}
The Hopf image of $A$ by a representation $\pi:A\to M_n(\mathbb C)$ is the smallest Hopf $C^*$-algebra $A'$ producing a factorization $\pi:A\to A'\to M_n(\mathbb C)$.
\end{definition}

In order to understand this notion, let $A=C^*(\Gamma)$. Then $\pi$ must come from a unitary group representation $\tilde{\pi}:\Gamma\to U_n$, and we have $A'=C^*(\Lambda)$, where $\Lambda=\tilde{\pi}(\Gamma)$.

In this computation $\Gamma$ was of course a usual discrete group. In the general case, i.e. when $\Gamma$ is a discrete quantum group, it is only known that the Hopf image exists, and is unique \cite{ifr}. But, of course, the discrete quantum group point of view is very useful.

Here are a few more definitions from \cite{ifr}, based on the same philosophy:

\begin{definition}
Let $A$ be a Hopf $C^*$-algebra.
\begin{enumerate}
\item A representation $\pi:A\to M_n(\mathbb C)$ is called inner faithful if $A=A'$.

\item $A$ is called inner linear if it has an inner faithful representation.
\end{enumerate}
\end{definition}

Observe that with $A=C^*(\Gamma)$, and with the above notations, $\pi$ is inner faithful if and only if $\tilde{\pi}$ is faithful. Also, $A$ is inner linear if and only if $\Gamma$ is linear.

These notions, coming from \cite{ver}, \cite{hdm}, \cite{roq}, were axiomatized and studied in \cite{ifr}. A number of key examples were later on studied by Andruskiewitsch and Bichon in \cite{abi}. 

Some useful supplementary tools, of functional analytic nature, are expected to come from the idempotent state theory of Franz and Skalski \cite{fsk}.

We can apply these constructions to the representation $\pi_H$ in Proposition 4.2:

\begin{definition}
Associated to any complex Hadamard matrix $H\in M_n(\mathbb C)$ is the quantum permutation group $G\subset S_n^+$ given by the fact that $C(G)$ is the Hopf image of $\pi_H$.
\end{definition}

In other words, if we go back to the definition of both $\pi_H$ and of the Hopf image, $C(G)$ is the smallest Hopf $C^*$-algebra producing the following factorization:
$$C(S_n^+)\to C(G)\to M_n(\mathbb C)$$
$$u_{ij}\to Proj(H_i/H_j)$$

Equivalently, we can say that this representation $\pi_H$ should correspond to a unitary quantum group representation $\tilde{\pi}_H:\widehat{S}_n^+\to U_n^+$, and $G$ is simply given by:
$$\widehat{G}=\tilde{\pi}_H(\widehat{S}_n^+)$$

This latter description of the correspondence $H\to G$, while being quite heuristic, starts to justify the first claim in Fact 5.1, namely that ``$G$ encodes the symmetries of $H$''.

\begin{proposition}
The construction $H\to G$ has the following properties:
\begin{enumerate}
\item For $H=F_n$ we obtain $G=\mathbb Z_n$.

\item For $H=H'\otimes H''$ we obtain $G=G'\times G''$.
\end{enumerate}
\end{proposition} 

The first assertion follows indeed from the definition of the Hopf image, because the  representation $C(S_n^+)\to M_n(\mathbb C)$ associated to $F_n$ is circulant, and hence factorizes through $C(\mathbb Z_n)$. Now since this factorization cannot be further factorized, we obtain $G=\mathbb Z_n$. 

As for the second assertion, this follows as well from definitions. See \cite{hdm}, \cite{roq}.

More generally, we have the following result, once again from \cite{hdm}, \cite{roq}:

\begin{theorem}
For an Hadamard matrix $H$, the following are equivalent:
\begin{enumerate}
\item $G$ is a classical group.

\item $\widehat{G}$ is a classical group.

\item $G=\mathbb Z_{n_1}\times\ldots\times\mathbb Z_{n_k}$, for some numbers $n_1,\ldots,n_k$.

\item $H=F_{n_1}\otimes\ldots\otimes F_{n_k}$, for some numbers $n_1,\ldots,n_k$.
\end{enumerate}
\end{theorem}

We recall now from Definition 4.5 above that the Di\c t\u a deformation of a tensor product $H\otimes K\in M_{nm}(\mathbb C)$, with matrix of parameters $L\in M_{m\times n}(\mathbb T)$, is:
$$H\otimes_LK=(H_{ij}L_{aj}K_{ab})_{ia,jb}$$

We have already seen in section 4 that this notion is a key one. For instance Theorem 4.4 tells us that at $n=4$ the only Hadamard matrices are the Di\c t\u a deformations of $F_2\otimes F_2$, and Theorem 4.8 tells us that the only regular matrices at $n=6$ are the Haagerup and Tao matrices $H^q,T$, and the Di\c t\u a deformations of $F_6=F_2\otimes F_3=F_3\otimes F_2$. 

The symmetries of $H\otimes_LK$, even for very simple values of $H,K$, vary a lot, and in a very subtle way, with the matrix of parameters $L$. The only general result here is:

\begin{theorem}
We have $G(H\otimes_LK)\subset G(K)\wr_*G(H)$.
\end{theorem}

This result, from \cite{roq}, marks the end of the algebraic theory. We will come back in section 6 below with a conjectural discussion of the problem.

For the moment, let us just record a difficult question:

\begin{problem}
What is $G(F_n\otimes_LF_m)$?
\end{problem}

Here are a few more questions, for the most of algebraic nature:

\begin{problem}
Let $H$ be a complex Hadamard matrix.
\begin{enumerate}
\item When is $G$ finite?

\item When is $\widehat{G}$ amenable?

\item Can we have $G=S_n^+$?
\end{enumerate}
\end{problem}

All these questions look equally very difficult. Question (3) is due to Jones \cite{jo3}.  Let us discuss now the representation theory of $G$. The result here, from \cite{roq}, is as follows:

\begin{theorem}
We have $T\in Hom(u^{\otimes k},u^{\otimes l})$ if and only if $T^\circ G^{k+2}=G^{l+2}T^\circ$, where:
\begin{enumerate}
\item $T^\circ=id\otimes T\otimes id$.

\item $G_{ia}^{jb}=\sum_{k=1}^nH_{ik}\bar{H}_{jk}\bar{H}_{ak}H_{bk}$.

\item $G^k_{i_1\ldots i_k,j_1\ldots j_k}=G_{i_ki_{k-1}}^{j_kj_{k-1}}\ldots G_{i_2i_1}^{j_2j_1}$.
\end{enumerate}
\end{theorem}

The proof in \cite{roq} uses a quite simple Tannakian argument, further developed in \cite{ifr}, stating that when we have an inner faithful representation $u_{ij}\to P_{ij}$, then:
$$Hom(u^{\otimes k},u^{\otimes l})=Hom(P^{\otimes k},P^{\otimes l})$$

Here the spaces on the right are defined exactly as those on the left, namely:
$$Hom(P^{\otimes k},P^{\otimes l})=\{T:(\mathbb C^n)^{\otimes k}\to(\mathbb C^n)^{\otimes l}|TP^{\otimes k}=P^{\otimes l}T\}$$

This result makes the link with the planar algebra approach of Jones \cite{jo3}:

\begin{theorem}
Let $H$ be a complex Hadamard matrix.
\begin{enumerate}
\item The planar algebra associated to $H$ is given by $P_k=Fix(u^{\otimes k})$.

\item The quantum invariants of $H$ are the moments of $\chi=Tr(u)$.

\item The Poincar\'e series of $H$ is the Stieltjes transform of $law(\chi)$.
\end{enumerate}
\end{theorem}

Indeed, (1) follows Theorem 5.10 at $k=0$, because the description there of the space $Fix(u^{\otimes l})$ coincides with the description in \cite{jo3} of the space $P_l$. As for  (2,3), these are simply quantum group translations of Jones' notions in \cite{jo3}, namely:
$$c_k=\dim(P_k),\quad\quad f(z)=\sum_{k=0}^\infty c_kz^k$$ 

More precisely, using the results of Woronowicz in \cite{wo1}, we get as claimed:
$$c_k=\int_G\chi^k,\quad\quad f(x)=\int_G\frac{1}{1-z\chi}$$ 

Summarizing, the above results fully justify the first claim in Fact 5.1, namely that ``$G$ is the symmetry group of $H$''. The story here can be actually made longer, by talking about commuting squares \cite{pop}, spin models \cite{jo2} and subfactors \cite{jo1}. See \cite{ver}.

We should probably mention that, via the various algebraic results presented in the beginning of this section, Theorem 5.11 allows one to recover most of what is known about the Hadamard matrix subfactors, sometimes with some slight enhancements. 

Let us go now to the second claim in Fact 5.1, namely that any quantum permutation group $G\subset S_n^+$ is the symmetry group of an Hadamard-type object $H\in\tilde{M}_n(\mathbb C)$. The situation here is much more complicated, and we only have a categorical result so far:

\begin{theorem}
There is a Tannaka-Galois correspondence between quantum permutation groups $G\subset S_n^+$ and subalgebras of the spin planar algebra $P\subset Q_n$. In addition:
\begin{enumerate}
\item The quantum symmetry groups of generalized colored graphs correspond in this way to the subalgebras of $Q_n$ generated by a $2$-box.

\item The quantum groups $S_n^+,H_n^+$ correspond in this way to the canonical copies of the Temperley-Lieb and Fuss-Catalan algebras inside $Q_n$.

\item The quantum group associated to a complex Hadamard matrix $H\in M_n(\mathbb C)$ corresponds in this way to the planar algebra associated to $H$. 
\end{enumerate}
\end{theorem}

The Tannaka-Galois correspondence was established in \cite{hgr}, the idea being that the embedding $P\subset Q_n$ could be regarded as a ``piece'' of the usual Tannakian data, namely ``tensor category + Hilbert functor''. By Frobenius duality we can recover the rest of the Tannakian data, and by using Woronowicz's duality in \cite{wo2} we obtain the result.

The assertion (1) was equally established in \cite{hgr}, the idea here being simply that the condition $d\in End(u)$ reads $d\in P_2$, and hence $P=<d>$. This result shows that our classification program for graphs, presented in section 3, can be regarded as part of the Bisch-Jones classification program for planar algebras generated by a 2-box \cite{bj2}.

The assertion (2) for $S_n^+$ follows simply from functoriality: since $S_n^+$ is the biggest quantum permutation group, it should correspond to the smallest planar subalgebra of $Q_n$, which is of course the Temperley-Lieb algebra $TL_n$. Observe that we have here yet another proof of Theorem 2.7, but this time, without any single computation! 

As for the assertion (2) for $H_n^+$, this follows from Theorem 3.12 above, or rather from the diagrammatic comments following it. See \cite{ahn}. The interest in this statement comes of course in conjunction with Theorem 3.12, and with the various probabilistic formulae following it: all these formulae concern of course the Fuss-Catalan algebra of Bisch and Jones \cite{bj1}, which was previously not known to have all these combinatorial aspects.

Finally, the assertion (3) follows from Theorem 5.11 above.

Let us go back now to the second claim in Fact 5.1. Given a quantum permutation group $G\subset S_n^+$, we would like to construct an object $H\in\tilde{M}_n(\mathbb C)$ whose symmetry group is $G$. This is of course a quite difficult task, and the above categorical result, or the purely algebraic considerations in \cite{roq}, are just a first step towards an answer.

We believe that the following question can push things here one step further:

\begin{problem}
What is the noncommutative geometry of a complex Hadamard matrix?
\end{problem}

In other words, associated to any complex Hadamard matrix $H\in M_n(\mathbb C)$ we would like to have a noncommutative manifold $M_H$ in the sense of Connes \cite{con}, modelled on the noncommutative set $M_n=Spec(M_n(\mathbb C))$, with the  requirement $G(H)=G^+(M)$.

\section{Matrix models}

We have seen in the previous section that associated to any complex Hadamard matrix $H\in M_n(\mathbb C)$ is a quantum permutation group $G\subset S_n^+$, describing its ``symmetries''. The correspondence $H\to G$ appears to be fully satisfactory at the level of simple examples, and is also compatible with the planar algebra approach of Jones \cite{jo3}.

The challenging problem is to relate the combinatorics of quantum permutation groups, explained in sections 1-3 above, to the combinatorics of complex Hadamard matrices, explained in section 4 above. Our answer here, largely conjectural, is as follows:

\begin{conjecture}
Consider the correspondence $H\to G$.
\begin{enumerate}
\item The known combinatorics of $G$ and the known combinatorics of $H$ are not different: they are just different aspects of a single combinatorics.

\item This single combinatorics appears when looking at the affine families $\{H^q|q\in\mathbb T^r\}$, and at the corresponding quantum groups $\{G^q|q\in\mathbb T^r\}$.

\item The unifying theorem is a statement regarding matrix models for $G^q,H^q$, with the parameter space being a certain subgroup of $\mathbb T^r$, associated to $q$.
\end{enumerate}
\end{conjecture}

In other words, our claim is that the quantum permutation groups $G^q=G(H^q)$ with $q\in\mathbb T^r$ should behave a bit like the usual Drinfeld-Jimbo deformations \cite{dri}, \cite{jim} in the parameter range $|q|=1$, with a rich combinatorics depending on the arithmetics of $q$. And this combinatorics should correspond to the various combinatorial aspects of $H^q$, i.e. to the Butson obstructions, and to the subtle behavior of the Poincar\'e series.

Our second claim is that a useful and effective tool for dealing with these difficult questions in the notion of matrix model. The idea here, to be explained in detail later on, is that under suitable assumptions we should have an explicit matrix model of type $C(G^q)\subset C(K_q)\otimes M_n(\mathbb C)$, where $K_q\subset\mathbb T^r$ is a certain abelian group associated to $q$.

Let us first explain the notion of affine family:

\begin{definition}
An affine family of complex Hadamard matrices is a family $H^q\in M_n(\mathbb C)$, depending analytically on a parameter $q\in\mathbb T^r$. The family is called:
\begin{enumerate}
\item Regular, if all matrices $H^q$ are regular.

\item Normalized, if $l(H^1)\leq l(H^q)$, for any $q\in\mathbb T^r$.

\item Arithmetic, if $l(H^q)\leq K\cdot ord(q)$, for some $K>0$.
\end{enumerate}
\end{definition}

We use here the notions of regularity and level, introduced in section 4. As an example, the Di\c t\u a deformations of Fourier matrices, as well as the Haagerup and Petrescu matrices $H^q,P^q$, form affine families which are regular, arithmetic, and often normalized. 

With a bit of complex analysis, one can prove that any arithmetic regular family is obtained in a similar way, by starting with a Butson matrix, and adding $r$ parameters.

We recall from \cite{nic}, \cite{tz2} that the ``defect'' of a matrix $H\in M_n(\mathbb C)$ counts the dimension of the component of the Hadamard matrix manifold $H_n(\infty)$ passing through $H$. 

Here is a first question that we have, regarding the affine families:

\begin{problem}
Is there a notion of ``regular defect'' of $H\in M_n(\mathbb C)$, counting the dimension of the maximal regular affine family passing through $H$?
\end{problem}

In other words, we are asking here for an explicit formula for the regular defect, in the spirit of those found in \cite{tz2}. This question is probably quite difficult, due to its relation with Conjecture 4.16 above, stating that any Butson matrix is regular. 

Here is another question, probably more elementary:

\begin{problem}
Is it true that the maximal affine family passing through $F=F_{n_1}\otimes\ldots\otimes F_{n_k}$ can be obtained via a sequence of Di\c t\u a deformations of $F$? 
\end{problem}

Normally this should follow from the various defect formulae in \cite{tz2} and from the fact that the Fourier matrix $F_{n_1\ldots n_k}$ is in the orbit of $F$, cf. Proposition 4.6.

Let us discuss now the notion of matrix model:

\begin{definition}
A matrix model for a Hopf $C^*$-algebra $A$ is a morphism of $C^*$-algebras $\pi:A\to L^\infty(X)\otimes B(H)$, where $X$ is a probability space, and $H$ is a finite dimensional Hilbert space, such that $\int a=(\mathbb E\otimes tr)\pi(a)$, for any $a\in A$.
\end{definition}

At the level of examples, let us first look at the commutative case. Here the algebras in question are those of the form $A=C(G)$, where $G$ is a usual compact group, and the identity map $id:C(G)\to L^\infty(G)\otimes B(\mathbb C)$ is of course a matrix model.

Let us record now a few general properties of the matrix models. 

\begin{proposition}
Let $\pi:A\to L^\infty(X)\otimes B(H)$ be a matrix model.
\begin{enumerate}
\item $A$ is amenable in the discrete quantum group sense.

\item $\pi$ is injective, i.e. $\pi(a)=0$ implies $a=0$.
\end{enumerate}
\end{proposition}

Indeed, let $\varphi:A\to A_{red}$ the projection onto the reduced version, obtained by dividing by the null ideal of the Haar functional. For any $a\in A$ we have:
$$\pi(a)=0
\implies(\mathbb E\otimes tr)\pi(aa^*)=0
\implies\int aa^*=0
\implies\varphi(a)=0$$

Thus we have $\ker\pi\subset\ker\varphi$, so $A_{red}$ appears as quotient of the image algebra $\pi(A)\subset L^\infty(X)\otimes B(H)$. With this observation in hand, we can prove now our assertions:

Indeed, (1) follows from the fact that $L^\infty(X)\otimes B(H)$ is of type I: by the general theory, so must be its subalgebra $\pi(A)$, and also its subquotient $A_{red}$. But now since $A_{red}$ is of type I, $A$ follows to be amenable, i.e. we have $A=A_{red}$, and we are done.

As for (2), since $A$ is amenable we have $\ker\pi\subset\ker\varphi=\mathbb C$, and we are done.

\begin{definition}
A linear model for a quantum permutation algebra $C(S_n^+)\to A$ is a continuous matrix model $\pi:A\to C(K)\otimes M_n(\mathbb C)$, where $K\subset O_n$ is a compact group, such that the magic entries $U_{ij}^x\in M_n(\mathbb C)$ are linear functions of $x\in K\subset M_n(\mathbb C)$.
\end{definition}

In other words, we are making some key assumptions here, both taking into account the number $n$ of points on which $A$ coacts: we assume that our probability space is a compact group $K\subset O_n$, that our Hilbert space is $H=\mathbb C^n$, and that the model itself comes from a ``linear construction'' applied to the embedding $K\subset B(H)$.

The basic example is the model for $C(S_4^+)$ found in \cite{iop}. Recall first that the elements of $SU_2$ are of the following form, with $(x,y,z,t)\in S^3$:
$$m=\begin{pmatrix}
x+iy&z+it\\
-z+it&x-iy
\end{pmatrix}$$

Observe that $m=xc_1+yc_2+zc_3+tc_4$, where $c_i$ are the Pauli matrices: 
$$c_1=\begin{pmatrix}1&0\cr 0&1\end{pmatrix}\hskip 5mm
c_2=\begin{pmatrix}i&0\cr 0&-i\end{pmatrix}\hskip 5mm
c_3=\begin{pmatrix}0&1\cr -1&0\end{pmatrix}\hskip 5mm
c_4=\begin{pmatrix}0&i\cr i&0\end{pmatrix}$$

Now $SU_2$ acts by conjugation, $m'(p)=mpm^*$, on $span(c_1,c_2,c_3,c_4)\simeq\mathbb R^4$, and we have $m'=diag(1,m'')$, with $m''\in SO_3$. Thus we have explicit morphisms $SU_2\to SO_3\subset O_4$.

\begin{theorem}
We have a linear model $C(S_4^+)\subset C(SO_3)\otimes M_4(\mathbb C)$ given by $U_{ij}^x=Proj(c_ixc_j)$, where $c_1,c_2,c_3,c_4$ are the Pauli matrices.
\end{theorem}

This result was found in \cite{iop}. The first observation is that the formula of $U_{ij}^x$ in the statement defines a morphism of algebras $C(S_4^+)\to C(SU_2)\otimes M_4(\mathbb C)$. Now since $U_{ij}^x$ is invariant under $x\to -x$, this shows that the target of our  morphism is contained in $C(SO_3)\otimes M_4(\mathbb C)$. Thus we have a representation as in the statement. 

As pointed out in \cite{qfp}, this representation is indeed a linear model in the sense of Definition 6.7, when using the above-mentioned embedding $SO_3\subset O_4$.

The proof in \cite{iop} of the equality $\int a=(\mathbb E\otimes tr)\pi(a)$ uses the Weingarten formula for the evaluation of the left term, and an explicit computation for the right term.

We have now all ingredients for stating our main conjecture:

\begin{conjecture}
Let $\{H^q|q\in\mathbb T^r\}$ be an affine family of complex Hadamard matrices, and consider the corresponding family of quantum permutation groups $\{G^q|q\in\mathbb T^r\}$. Then, under suitable assumptions, there exist compact abelian groups $K_q\subset\mathbb T^r$ such that:
\begin{enumerate}
\item We have twisting results of type $G_q=K_q^{-1}$.

\item We have linear models of type $C(G_q)\subset C(K_q)\otimes M_n(\mathbb C)$.
\end{enumerate}
\end{conjecture}

Observe that at $n=4$ we know from Theorem 5.7 above that we have an inclusion $G(F_4^q)\subset\mathbb Z_2\wr_*\mathbb Z_2=O_2^{-1}$, so the remaining problem is to find the correct group $K_q$, by using Theorem 1.12 and Theorem 6.8. This will be explained in a forthcoming paper.

\section*{Conclusion}

We have seen in this paper that the combinatorics of quantum permutation groups $G\subset S_n^+$ shares many common aspects with the combinatorics of complex Hadamard matrices $H\in M_n(\mathbb C)$. The quite substantial mass of work so far, on quantum groups and Hadamard matrices, goes however is slightly different directions. We can only hope that a further study of the correspondence $H\to G$ would be fruitful on both sides.

One problem comes from the fact that a lot of work on $G$ has gone into quite delicate probabilistic aspects, whose matrix counterpart is for the moment unknown.

The other problem is that a lot of work on $H$ has gone into the root of unity case, and an abstract notion of quantum permutation groups at roots of unity is still lacking.

We believe that both problems can be solved, and the correspondence $H\to G$ can be substantially upgraded, in the general context of the matrix model problematics. This idea, which is largely conjectural for the moment, is explained in section 6 above.

Finally, we hope that the present text was readable. The prequisites for a second reading would be a certain familiarity with Hopf algebras and operator algebras \cite{abe}, \cite{ped}, with compact quantum groups \cite{wo1}, \cite{wo2}, and with free probability theory \cite{nsp}, \cite{vdn}. The reader might also consult the previous survey article on the subject \cite{qpg}.

\end{document}